\documentstyle[12pt]{article}  
\def\sq{\hbox {\rlap{$\sqcap$}$\sqcup$}}
\def\sq{\hbox {\rlap{$\sqcap$}$\sqcup$}}
\def\R{ {\rm R \kern -.31cm I \kern .15cm}}
\def\C{ {\rm C \kern -.15cm \vrule width.5pt \kern .12cm}}
\def\Z{ {\rm Z \kern -.27cm \angle \kern .02cm}}
\def\N{ {\rm N \kern -.26cm \vrule width.4pt \kern .10cm}}
\def\1{{\rm 1\mskip-4.5mu l} }
\def\lsim{\raise0.3ex\hbox{$<$\kern-0.75em\raise-1.1ex\hbox{$\sim$}}}
\def\gsim{\raise0.3ex\hbox{$>$\kern-0.75em\raise-1.1ex\hbox{$\sim$}}}
\def\noi{\noindent}

\def\beq{\begin{equation}}   \def\eeq{\end{equation}}
\def\bea{\begin{eqnarray}}  \def\eea{\end{eqnarray}}
\def\nn{\nonumber}
\def\noi{\noindent}
\def\beeq{\begin{eqnarray}} \def\eeeq{\end{eqnarray}}
\newcommand\mysection{\setcounter{equation}{0}\section}

\newcounter{hran}

\begin{document} 
\centerline{\large\bf Scattering theory for the Zakharov system} 
\vskip 0.5 truecm

\centerline{\bf J. Ginibre}
\centerline{Laboratoire de Physique Th\'eorique\footnote{Unit\'e Mixte de
Recherche (CNRS) UMR 8627}}  \centerline{Universit\'e de Paris XI, B\^atiment
210, F-91405 ORSAY Cedex, France}
\vskip 3 truemm
\centerline{\bf G. Velo}
\centerline{Dipartimento di Fisica, Universit\`a di Bologna}  \centerline{and INFN, Sezione di
Bologna, Italy}

\vskip 1 truecm

\begin{abstract}
We study the theory of scattering for the Zakharov system in space
dimension 3. We prove in particular the existence of wave
operators for that system with no size restriction on the data in
larger spaces and for more general asymptotic states than were
previously considered, and we determine convergence rates  in
time of solutions in the range of the wave operators to the solutions of the underlying linear system. We also consider
the same system in space dimension 2, where we prove the existence of
wave operators for small Schr\"odinger data in the special case of
vanishing asymptotic data for the wave field.
\end{abstract}

\vskip 3 truecm
\noi AMS Classification : Primary 35P25. Secondary 35B40, 35Q60, 81U99.  \par \vskip 2 truemm

\noi Key words : Scattering theory, Zakharov system.\par 
\vskip 1 truecm

\noindent LPT Orsay 05-19\par
\noindent March 2005\par

\newpage
\pagestyle{plain}
\baselineskip 18pt

\mysection{Introduction}
\hspace*{\parindent} This paper is devoted to the theory of scattering
and more precisely to the construction of wave operators for
the Zakharov system (Z)$_n$ in space dimension $n = 3$ and 2 (in that order), namely
\beq \label{1.1e}
\left \{ \begin{array}{l}  i\partial_t u = - (1/2)
\Delta u + A u \\ \\ \sq A = \Delta |u|^2   \end{array} \right . 
 \eeq

\noi where $u$ and $A$ are respectively a complex valued and a real
valued function defined in space time ${I\hskip-1truemm R}^{n+1}$,
$\Delta$ is the Laplacian in ${I\hskip-1truemm R}^n$ and $\sq =
\partial_t^2 - \Delta$ is the d'Alembertian in ${I\hskip-1truemm
R}^{n+1}$. The $(Z)_3$ system is used in plasma physics to describe the
Langmuir turbulence. The function $u$ is the slowly varying complex
amplitude of the rapidly oscillating electric field and $A$ is the
deviation of the ion density from its average value \cite{27r}. In this
paper we  use the notation $(u, A)$ for those variables instead of the
more common $(E, n)$ in order to allow for a better contact with the
existing literature on related nonlinear systems based on the
Schr\"odinger equation, in particular with the Wave-Schr\"odinger
system (WS)$_3$ and the Maxwell-Schr\"odinger system (MS)$_3$ in
${I\hskip-1truemm R}^{3+1}$, and with the Klein-Gordon-Schr\"odinger
system (KGS)$_2$ in ${I\hskip-1truemm R}^{2+1}$ (see \cite{13r} for a
review). The Zakharov system is Lagrangian and admits a number of
formally conserved quantities, among which the $L^2$ norm of $u$ and
the energy
\beq
\label{1.2e}
E(u,A) = \int dx \left \{ (1/2) \left ( |\nabla u|^2 + |\omega^{-1}\partial_t A|^2 + |A|^2\right ) + A|u|^2 \right \} 
\eeq

\noi where $\omega = (- \Delta )^{1/2}$. The Cauchy problem for the
(Z)$_n$ system has been extensively studied \cite{1r} \cite{2r}
\cite{3r} \cite{6r} \cite{15r} \cite{25r} and is known to be globally
well posed for $n = 2,3$ in the energy space $X_{e} = H^1 \oplus L^2
\oplus \dot{H}^{-1}$ for $(u, A, \partial_t A)$.\par

A large amount of work has been devoted to the theory of scattering for
non linear equations and systems related to the Schr\"odinger equation,
in particular for non linear Schr\"odinger (NLS)$_n$ and Hartree
(R3)$_n$ equations in ${I\hskip-1truemm R}^{n+1}$ and for the above
mentioned (WS)$_3$, (MS)$_3$, (KGS)$_2$ and (Z)$_3$ systems. As in the
case of the linear Schr\"odinger equation, one must distinguish the
short range case from the long range case. In the former case, ordinary
wave operators are expected and in a number of cases proved to exist,
describing solutions where the Schr\"odinger function behaves
asymptotically like a solution of the free Schr\"odinger equation. In
the latter case, ordinary wave operators do not exist and have to be
replaced by modified wave operators including a suitable phase in their
definition. In that respect, the (WS)$_3$ and (MS)$_3$ systems belong
to the borderline (Coulomb) long range case, as does the (R3)$_n$
equation with $|x|^{-1}$ potential, the (Z)$_3$ system is short range,
and the (KGS)$_2$ and (Z)$_2$ systems, although not really long range,
exhibit some difficulties typical of the long range case. \par

The construction of (possibly modified) wave operators for the
previous equations and systems in the long range cases has been tackled
by two methods. The first one was initiated on the example of (NLS)$_1$
\cite{14r} and subsequently applied to the (NLS)$_n$ equation for $n =
2,3$ and to the (R3)$_n$ equation for $n \geq 2$ \cite{5r}, to the
(KGS)$_2$ system \cite{17r} \cite{18r} \cite{19r} \cite{20r}, to the
(WS)$_3$ system \cite{11r} \cite{21r}, to the (MS)$_3$ system
\cite{12r} \cite{22r} \cite{26r} and to the (Z)$_3$ system \cite{16r}
\cite{23r}. See \cite{13r} for a review. That method is rather direct,
starting from the original equation or system. It will be sketched
below on the example of the (Z)$_n$ system. In long range cases, it is
restricted to the limiting Coulomb case and requires a smallness
condition on the asymptotic state of the Schr\"odinger function. Early
applications of the method required in addition a support condition on
the Fourier transform of the Schr\"odinger asymptotic state and a
smallness condition of the Klein-Gordon or Maxwell field in the case of
the (KGS)$_2$ or (MS)$_3$ system respectively \cite{17r} \cite{26r}. A
support condition was also required in the case of the (Z)$_3$ system
when both the Schr\"odinger and the wave field are large, which is
allowed by the fact that the (Z)$_3$ system is short range \cite{16r}.
The support condition was subsequently removed for the (KGS)$_2$,
(MS)$_3$ and (Z)$_3$ systems and the method was applied to the (WS)$_3$
system without a support condition, at the expense of adding a
correction term to the Schr\"odinger asymptotic function \cite{18r}
\cite{21r} \cite{22r} \cite{23r}. The smallness condition of
the KG field was then removed for the (KGS)$_2$ system, first with and
then without a support condition \cite{19r} \cite{20r}. All the
previous papers on (KGS)$_2$, (WS)$_3$, (MS)$_3$ and (Z)$_3$ use spaces
of fairly regular solutions, with at least $H^2$ regularity for the
Schr\"odinger function. Finally the smallness condition of the wave or
Maxwell field was removed for the (WS)$_3$ and (MS)$_3$ systems
\cite{11r} \cite{12r}. Furthermore larger function spaces than
previously considered are used in \cite{11r} \cite{12r}, thereby
allowing for more general asymptotic states.\par

In the present paper, we reconsider the same problem for the (Z)$_3$
and (Z)$_2$ systems in the framework of the previous method. We treat
again the (Z)$_3$ system with no smallness condition on either field
and no support condition. In the same spirit as in \cite{11r}
\cite{12r}, we use function spaces that are as large as possible, namely with regularity
as low as possible, and with convergence in time as slow as possible.
In particular we treat the problem with only a weak convergence in time
of the solutions to their asymptotic form, namely $t^{-\lambda}$ with
$\lambda > 1/4$. Under such a weak condition, neither a support
condition nor a correction term for the asymptotic Schr\"odinger
function is needed as long as $\lambda \leq 1/2$ and much weaker
assumptions on the asymptotic state than previously considered can be
accomodated. We also consider the case of more regular data but still
more general than previously considered, where the use of a correction
term yields a stronger convergence in time, namely $\lambda = 3/2$. We
finally apply the method to the (Z)$_2$ system. Again no support
condition is needed, but we need a smallness condition of the
Schr\"odinger function and we can only treat the case where the
asymptotic state of the wave field is zero.\par

For completeness and although we shall not make use of that fact in the
present paper, we mention that the same problem for the Hartree
equation and for the (WS)$_3$ and (MS)$_3$ system can also be treated
by a more complex method where one first applies a phase-amplitude
separation to the Schr\"odinger function. The main interest of that
method is to remove the smallness condition on the Schr\"odinger
function, and to go beyond the Coulomb limiting case for the Hartree
equation. That method has been applied in particular to the (WS)$_3$ system and to
the (MS)$_3$ system in a special case \cite{7r} \cite{8r} \cite{9r}.
\par

We now sketch briefly the method of construction of the modified wave
operators initiated in \cite{14r}. That construction basically consists
in solving the Cauchy problem for the system (\ref{1.1e}) with infinite
initial time, namely in constructing solutions $(u,A)$ with prescribed
asymptotic behaviour at infinity in time. We restrict our attention to
time going to $+\infty$. That asymptotic behaviour is imposed in the
form of suitable approximate solutions $(u_a,A_a)$ of the system
(\ref{1.1e}). The approximate solutions are parametrized by data $(u_+,
A_+, \dot{A}_+)$ which in the simplest cases are initial data at time zero for a simpler
evolution. One then looks for exact solutions $(u,A)$ of the system
(\ref{1.1e}), the difference of which with the given asymptotic ones
tends to zero at infinity in time in a suitable sense, more precisely,
in suitable norms. The wave operator is then defined traditionally as
the map $\Omega_+ : (u_+, A_+, \dot{A}_+) \to (u,A,\partial_t A)(0)$.
However what really matters is the solution $(u, A)$ in the
neighborhood of infinity in time, namely in some interval $[T, \infty
)$, and we shall restrict our attention to the construction of such
solutions. Continuing such solutions down to $t = 0$ is a somewhat
different question, connected with the global Cauchy problem at finite
times, which we shall not touch here. That problem is well controlled
for the (Z)$_n$ system for $n = 2,3$.
\par

The construction of solutions $(u,A)$ with prescribed asymptotic
behaviour $(u_a,A_a)$ is performed in two steps. \\

\noi \underbar{Step 1}. One looks for $(u,A)$ in the form $(u,A) = (u_a+v, A_a + B)$. The system satisfied by $(v,B)$ is
\beq \label{1.3e}
\left \{ \begin{array}{l}  i\partial_t v = - (1/2)
\Delta v + A v + Bu_a - R_1 \\ \\ \sq B = \Delta (|v|^2 + 2 \ {\rm Re} \ \overline{u}_a v) - R_2   \end{array} \right . 
 \eeq

\noi where the remainders $R_1$, $R_2$ are defined by
\beq \label{1.4e}
\left \{ \begin{array}{l}  R_1= i\partial_t u_a + (1/2)
\Delta u_a - A_a u_a \\ \\ R_2 = \sq A_a - \Delta |u_a|^2   \ . \end{array} \right . 
 \eeq

\noi It is technically useful to consider also the partly linearized system for functions $(v',B')$
\beq \label{1.5e}
\left \{ \begin{array}{l}  i\partial_t v' = - (1/2)
\Delta v' + A v' + Bu_a - R_1 \\ \\ \sq B' = \Delta (|v|^2 + 2 \ {\rm Re} \ \overline{u}_a v) - R_2  \ . \end{array} \right . 
 \eeq

\noi The first step of the method consists in solving the system
(\ref{1.3e}) for $(v, B)$, with $(v, B)$ tending to zero at infinity in
time in suitable norms, under assumptions on $(u_a, A_a)$ of a general
nature, the most important of which being decay assumptions on the
remainders $R_1$ and $R_2$. That can be done as follows. One first
solves the linearized system (\ref{1.5e}) for $(v',B')$ with given $(v,
B)$ and initial data $(v', B')(t_0) = 0$ for some large finite $t_0$.
One then takes the limit $t_0 \to \infty$ of that solution, thereby
obtaining a solution $(v',B')$ of (\ref{1.5e}) which tends to zero at
infinity in time. That construction defines a map $\phi : (v,B) \to
(v',B')$. One then shows by a contraction method that the map $\phi$
has a fixed point.  \\

\noi \underbar{Step 2.} The second step of the method consists in
constructing approximate asymptotic solutions $(u_a, A_a)$ satisfying
the general estimates needed to perform Step~1. With the weak time
decay allowed by our treatment of Step 1, and taking advantage of the
fact that the (Z)$_3$ system is short range, one can take for $(u_a,
A_a)$ solutions of the free Schr\"odinger and wave equations in that
case. One can also improve $u_a$ by a correction term as in \cite{23r},
thereby obtaining faster convergence rates for more regular asymptotic
states. In the case of the (Z)$_2$ system, one can again take for $u_a$
a solution of the free Schr\"odinger equation, but one is forced to
take $A_a = 0$.\par

In order to state our results we introduce some notation. We denote by
$F$ the Fourier transform and by $\parallel \cdot \parallel_r$ the norm
in $L^r \equiv L^r ({I\hskip-1truemm R}^n)$, $1 \leq r \leq \infty$.
For any nonnegative integer $k$ and for $1 \leq r \leq \infty$, we
denote by $W_r^k$ the Sobolev spaces 
$$W_r^k = \left \{ u : \parallel u; W_r^k\parallel \ = \sum_{\alpha : 0 \leq |\alpha | \leq k} \parallel \partial_x^{\alpha} u \parallel_r \ < \infty \right \}$$

\noi where $\alpha$ is a multiindex, so that $H^k = W_2^k$. We shall
need the weighted Sobolev spaces $H^{k,s}$ defined for $k$, $s \in
{I\hskip-1truemm R}$ by 
$$H^{k,s} = \left \{ u : \parallel u; H^{k,s}\parallel \ = \ \parallel (1 + x^2)^{s/2} (1 - \Delta )^{k/2} u \parallel_2 \ < \infty \right \}$$

\noi so that $H^k = H^{k,0}$. For any interval $I$, for any Banach
space $X$ and for any $q$, $1 \leq q \leq \infty$, we denote by $L^q(I,
X)$ (resp. $L_{loc}^q(I,X)$) the space of $L^q$ integrable (resp.
locally $L^q$ integrable) functions from $I$ to $X$ if $q < \infty$ and
the space of measurable essentially bounded (resp. locally essentially
bounded) functions from $I$ to $X$ if $q = \infty$. We shall occasionally use the notation 
$$\parallel f;L^q (I, L^r)\parallel\ = \ \parallel \ \parallel f \parallel_r \ \parallel_q$$

\noi when there is no ambiguity in the choice of the interval I. For any $h \in
{\cal C} ([1, \infty ), {I\hskip-1truemm R}^+)$, non increasing and
tending to zero at infinity and for any interval $I \subset [1, \infty
)$, we define the space
\bea
\label{1.6e}
&&X(I) = \Big \{ (v, B):(v,B)\in {\cal C}(I, H^2\oplus H^1) \cap {\cal C}^1(I, L^2 \oplus L^2), \nn \\
&&\parallel (v, B);X(I)\parallel \ \equiv \ \mathrel{\mathop {\rm Sup}_{t \in I }}\ h(t)^{-1} 
\Big ( \parallel v(t);H^2\parallel \ + \ \parallel \partial_t v(t)\parallel_2\nn \\
&&+\ \parallel v; L^{8/n}(J,W_4^2)\parallel \ + \ \parallel \partial_t v; L^{8/n}(J, L^4)\parallel \nn\\
&& + \ \parallel B(t); H^1\parallel\ + \ \parallel \partial_t B(t) \parallel_2 \Big ) < \infty \Big \} 
\eea

\noi where $J = [t, \infty ) \cap I$, for $n = 2,3$. Finally we denote by
\beq
\label{1.7e}
u_0(t) = U(t) u_+ = \exp (i(t/2)\Delta ) u_+ \ ,
\eeq
\beq
\label{1.8e}
A_0(t) = \cos \omega t \ A_+ + \omega^{-1} \sin \omega t \ \dot{A}_+
\eeq

\noi the solutions of the free Schr\"odinger and wave equations with
initial data $u_+$ and $(A_+, \dot{A}_+)$ at time zero. \par

We can now state our results. We first state the result obtained for
the (Z)$_3$ system by using only the simplest asymptotics (\ref{1.7e})
(\ref{1.8e}).\\

\noi {\bf Proposition 1.1.} {\it Let $n = 3$. Let $h(t) = t^{-1/2}$ and let $X(\cdot
)$ be defined by (\ref{1.6e}). Let $u_+ \in H^2 \cap W_1^2$, let $A_+,
\omega^{-1} \dot{A}_+ \in H^1$ and $\nabla^2 A_+, \nabla \dot{A}_+ \in
W_1^1$. Let $(u_0, A_0)$ be defined by (\ref{1.7e}) (\ref{1.8e}).
Then there exists $T$, $1 \leq T < \infty$, and there exists a unique
solution $(u,A)$ of the (Z)$_3$ system (\ref{1.1e}) such that $(v,B) \equiv
(u-u_0, A- A_0) \in X([T, \infty))$. If in addition $u_+ \in H^{0,2}$,
then $B$ satisfies the estimate
\beq
\label{1.9e}
\parallel B(t);H^1\parallel\ \vee \ \parallel \omega^{-1}\partial_t B(t); H^1)\parallel \ \leq C\ t^{-3/4}  
\eeq
 
 \noi for some constant $C$ and for all $t \geq T$.}\\
 
 We next state the result obtained for the (Z)$_3$ system by using an
improved asymptotic $u_a$, for more regular asymptotic states $(u_+,
A_+, \dot{A}_+)$ and with stronger asymptotic convergence in time.\\

\noi {\bf Proposition 1.2.} {\it Let $n=3$. Let $h(t) = t^{-3/2}$ and let $X(\cdot
)$ be defined by (\ref{1.6e}). Let $u_+ \in H^2 \cap H^{0,2} \cap
W_1^2$ with $xu_+ \in W_1^2$. Let $(A_+,\dot{A}_+)$ satisfy
\bea
\label{1.10e}
&&A_+, \omega^{-1} \dot{A}_+ \in \dot{H}^{-2} \cap H^1\quad , \quad \nabla^2A_+, \nabla \dot{A}_+ \in W_1^1 \ ,\nn\\
&&x\cdot \nabla A_+, \omega^{-1} x \cdot \nabla \dot{A}_+ \in \dot{H}^{-2} \cap \dot{H}^{-1}\ .
\eea

\noi Let $(u_0, A_0)$ be defined by (\ref{1.7e}) (\ref{1.8e}) and let
$u_a = (1 + f)u_0$ with $f = 2 \Delta^{-1}A_0$. Then~:\par

(1) There exists $T$,
$1 \leq T < \infty$ and there exists a unique solution $(u, A)$ of the
(Z)$_3$ system (\ref{1.1e}) such that $(v, B) \equiv (u - u_a, A - A_0) \in
X([T, \infty ))$.\par

(2) Assume in addition that $\omega^{-1} A_+$, $\omega^{-2} \dot{A}_+ \in W_{4/3}^1$. Then there exists $T$,
$1 \leq T < \infty$ and there exists a unique solution $(u, A)$ of the
(Z)$_3$ system (\ref{1.1e}) such that $(u - u_0, A - A_0) \in
X([T, \infty ))$. One can take the same $T$ and the solution $(u, A)$ is the same as in Part (1).}\\

We finally state the result for the (Z)$_2$ system. As already
mentioned, that result requires small Schr\"odinger data, namely small
$u_+$, and requires $A_+ = \dot{A}_+ = 0$.\\

\noi {\bf Proposition 1.3.} {\it Let $n = 2$. Let $h(t) = t^{-1}$ and let $X(\cdot
)$ be defined by (\ref{1.6e}). Let $u_+ \in H^2 \cap H^{0,2} \cap
W_1^2$ with $\parallel u_+ ; W_1^2\parallel$ sufficiently small and let $u_0(t) = U(t) u_+$.
 Then there exists $T$,
$1 \leq T < \infty$, and there exists a unique solution $(u, A)$ of the
(Z)$_2$ system (\ref{1.1e}) such that $(u - u_0, A ) \in
X([T, \infty ))$.}\\

\noi {\bf Remark 1.1.} We could have included the norm $\parallel
\omega^{-1} \partial_t B \parallel_2$, which is part of the energy, in
the definition of $X(\cdot )$. That norm is never used in the proofs to
perform the estimates and comes out at the end as a by product thereof.
We have omitted it for simplicity.\\

The results of this paper have been announced in \cite{13r}. 

\mysection{The Zakharov system (Z)$_{\bf 3}$ in space dimension n $=$ 3}
\hspace*{\parindent} In this section we treat the (Z)$_3$ system and
eventually prove Propositions 1.1 and 1.2. We follow the sketch given
in the introduction and begin with the first step of the method. The
treatment of that step follows exactly the same pattern as for the
(WS)$_3$ system treated in \cite{11r}. We shall therefore be rather
sketchy as regards the general arguments of the proofs, for which we
refer to \cite{11r} for more details, and we shall mostly concentrate
on the parts which are specific to the (Z)$_3$ system, namely the
estimates. We shall make extensive use of the Strichartz inequalities
for the Schr\"odinger equation \cite{4r} which we recall for
completeness, in space dimension $n \geq 2$. A pair of exponents $q$,
$r$ with $2 \leq q$, $r \leq \infty$ is called admissible if  
\beq
\label{2.1e}
\begin{array}{ll} 0 \leq 2/q = n/2 - n/r &\leq 1 \quad \hbox{for $n \geq 3$}\\
&\\
&<1 \quad \hbox{for $n = 2$} \ . \\
\end{array}
\eeq

\noi {\bf Lemma 2.1.} {\it Let $(q_i, r_i)$, $i = 1,2$, be two admissible pairs. Let $v$ satisfy the equation
$$i \partial_t v = - (1/2) \Delta v + f$$

\noi in some interval $I$ with $v(t_0) = v_0$ for some $t_0 \in I$. Then the following estimates hold~:
\beq
\label{2.2e}
\parallel v; L^{q_1}(I, L^{r_1})\parallel \ \leq \ C \left ( \parallel v_0\parallel_2 \ + \ \parallel f; L^{\overline{q}_2}(I, L^{\overline{r}_2})\parallel \right )
\eeq

\noi where $C$ is a constant independent of $I$, and with $1/p +
1/\overline{p} = 1$. }\\

Note that the pair $(8/n, 4)$ which appears in the definition (\ref{1.6e}) of $X(\cdot )$ is an admissible pair.\par

We shall also need some information on the Cauchy problem at finite
times for the Schr\"odinger equation with time dependent real potential
and time dependent inhomogeneity.
\beq
\label{2.3e}
i \partial_t v = - (1/2) \Delta v + Vv + f \ .
\eeq

\noi We refer to Proposition 3.2 in \cite{10r} for sufficient
conditions on $V$, $f$ under which that problem is (globally) well
posed with solutions in ${\cal C}(\cdot , H^2) \cap {\cal C}^1(\cdot ,
L^2)$.\par

We shall make repeated use of the following lemma, which is proved in \cite{11r}. \\

\noi {\bf Lemma 2.2.} {\it Let $1 \leq T < t_0 \leq \infty$, let $I = [T, t_0)$, and for $t \in I$, let $J = [t, \infty ) \cap I$. Let $1 \leq q$, $q_k \leq \infty$ ($1 \leq k \leq n$) be such that 
$$\mu \equiv 1/q - \sum_k 1/q_k \geq 0 \ .$$

\noi Let $f_k \in L^{q_k}(I)$ satisfy
\beq
\label{2.4e}
\parallel f_k;L^{q_k}(J) \parallel\ \leq N_k\ h(t)
\eeq

\noi for $1 \leq k \leq n$, for some constants $N_k$ and for all $t \in I$. \par

Let $\rho \geq 0$ be such that $n\lambda + \rho > \mu$. Then the following inequality holds for all $t\in I$
\beq
\label{2.5e}
\parallel \left ( \prod_k f_k \right ) t^{-\rho} ; L^q(J) \parallel\ \leq C \left ( \prod_k N_k\right ) h(t)^n \ t^{\mu - \rho}
\eeq

\noi where}
$$C = \left ( 1 - 2^{-q(n\lambda + \rho - \mu )}\right )^{-1/q} \ .$$
\vskip 5 truemm

We now turn to Step 1 of the method, namely to the construction of
solutions of the system (\ref{1.3e}) under general assumptions on
$(u_a, A_a)$. We denote by $h$ a function in ${\cal C}([1, \infty ), {I\hskip-1truemm R}^{+})$
such that for some $\lambda > 0$, the function $\overline{h}(t) =
t^{\lambda} h(t)$ is nonincreasing and tends to zero at infinity. The
main result on Step 1 can be stated as follows.\\

\noi {\bf Proposition 2.1} {\it Let $h$ be defined as above with
$\lambda = 1/4$, and let $X(\cdot )$ be defined by (\ref{1.6e}). Let
$u_a$, $A_a$, $R_1$ and $R_2$ be sufficiently regular (for the
following estimates to make sense) and satisfy the estimates
\beq
\label{2.6e}
\parallel u_a(t) \parallel_{\infty}\ \vee \ \parallel \nabla u_a(t) \parallel_{\infty} \ \vee \ \parallel \Delta u_a(t)\parallel_{\infty}\ \vee \ \parallel \partial_t u_a(t) \parallel_{\infty} \ \leq c\ t^{-3/2} 
\eeq
\beq
\label{2.7e}
\parallel \partial_t^j A_a(t) \parallel_{\infty} \ \leq a \ t^{-1} \qquad \hbox{\it for $j = 0, 1$}\ ,
\eeq  
\beq
\label{2.8e}
\parallel \partial_t^j R_1; L^1([t , \infty ), L^2)\parallel \ \leq r_1 \ h(t) \qquad \hbox{\it for $j = 0, 1$}\ ,
\eeq 
\beq
\label{2.9e}
\parallel R_1; L^{8/3}([t , \infty ), L^4)\parallel \ \leq r_1 \ t^{-\eta}\ h(t) \qquad \hbox{\it for some $\eta \geq 0$}\ ,
\eeq 
\beq
\label{2.10e}
\parallel \omega^{-1}R_2; L^1([t , \infty ), H^1)\parallel \ \leq r_2 \ h(t) \ ,
\eeq

\noi for some constants $c$, $a$, $r_1$ and $r_2$ and for all $t \geq
1$. Then there exists $T$, $1 \leq T < \infty$, and there exists a
unique solution $(v, B)$ of the system (\ref{1.3e}) in $X([T, \infty
))$. If in addition
\beq
\label{2.11e}
\parallel \omega^{-1}R_2; L^1([t , \infty ), L^2)\parallel \ \leq r_2 \ t^{-1/2}\ h(t) 
\eeq

\noi for all $t \geq T$, then $B$ satisfies the estimate
\beq
\label{2.12e}
\parallel B(t);H^1\parallel\ \vee \ \parallel \omega^{-1}\partial_t B(t); H^1)\parallel \ \leq C\left ( t^{-1/2} + t^{1/4}  \ h(t)\right ) \ h(t) 
\eeq

\noi for some constant $C$ and for all $t \geq T$.}\\

\noi {\bf Proof.} We follow the sketch given in the introduction, and more precisely the proof of Proposition 2.2 in \cite{11r}. Let $1 \leq T < \infty$
and let $(v, B) \in X([T, \infty))$. In particular $(v, B)$ satisfies the estimates 
\bea
\label{2.13e}
&&\parallel v(t) \parallel_2 \ \leq N_0\ h(t)\\
&&\parallel v ; L^{8/3}(J, L^4)\parallel \ \leq \  N_1\ h(t)\\ 
\label{2.14e}
&&\parallel B(t);H^1 \parallel \ \vee \ \parallel \partial_t B(t)\parallel_2\ \leq \  N_2\ h(t) \\
\label{2.15e}
&&\parallel \partial_t v(t) \parallel_2 \ \leq \  N_3\ h(t)\\
\label{2.16e}
&&\parallel \partial_t v ; L^{8/3}(J, L^4)\parallel \ \leq \  N_4\ h(t)\\ 
\label{2.17e}
&&\parallel \Delta v(t) \parallel_2\ \leq N_5\ h(t) \\
\label{2.18e}
&&\parallel \Delta v; L^{8/3}(J, L^4) \parallel\ \leq N_6\ h(t) 
\label{2.19e}
\eea

\noi for some constants $N_i$, $0 \leq i \leq 6$ and for all $t \geq
T$, with $J = [t, \infty )$. We first construct a solution $(v', B')$
of the system (\ref{1.5e}) in $X([T, \infty ))$. For that purpose, we
take $t_0$, $T < t_0 < \infty$ and we solve the system (\ref{1.5e}) in
$X(I)$ where $I = [T, t_0]$ with initial condition $(v', B')(t_0) = 0$.
Let $(v'_{t_0}, B'_{t_0})$ be the solution thereby obtained. The
existence of $v'_{t_0}$ follows from Proposition 3.2 in \cite{10r} with
$(A, V, f)$ replaced by $(0, A, - R_1)$. We want to take the limit of
$(v'_{t_0}, B'_{t_0})$ as $t_0 \to \infty$ and for that purpose we need
estimates of $(v'_{t_0}, B'_{t_0})$ in $X(I)$ that are uniform in
$t_0$. Omitting the subscript $t_0$ for brevity, we define
\beq
\label{2.20e}
N'_0 =  \ \mathrel{\mathop {\rm Sup}_{t\in I}}\ h(t)^{-1} \parallel v'(t) \parallel_2
\eeq 
\beq
\label{2.21e}
N'_1 =  \ \mathrel{\mathop {\rm Sup}_{t\in I}}\ h(t)^{-1} \parallel v';L^{8/3}(J,L^4) \parallel
\eeq
\beq
\label{2.22e}
N'_2 =  \ \mathrel{\mathop {\rm Sup}_{t\in I}}\ h(t)^{-1} \left ( \parallel B'(t);H^1 \parallel \ \vee \ \parallel \partial_t B'(t) \parallel_2 \right )
\eeq
\beq
\label{2.23e}
N'_3 =  \ \mathrel{\mathop {\rm Sup}_{t\in I}}\ h(t)^{-1} \parallel \partial_t v'(t)\parallel_2
\eeq
\beq
\label{2.24e}
N'_4 =  \ \mathrel{\mathop {\rm Sup}_{t\in I}}\ h(t)^{-1} \parallel \partial_t v';L^{8/3}(J,L^4) \parallel
\eeq
\beq
\label{2.25e}
N'_5 =  \ \mathrel{\mathop {\rm Sup}_{t\in I}}\ h(t)^{-1} \parallel \Delta v'(t)\parallel_2
\eeq
\beq
\label{2.26e}
N'_6 =  \ \mathrel{\mathop {\rm Sup}_{t\in I}}\ h(t)^{-1} \parallel \Delta v';L^{8/3}(J,L^4) \parallel
\eeq

\noi where $J = [t, \infty ) \cap I$ and we set out to estimate the various $N'_i$. We first estimate $N'_0$. From (\ref{1.5e}) we obtain
\bea
\label{2.27e}
\parallel v'(t)\parallel_2 &\leq& \parallel Bu_a - R_1;L^1(J, L^2)\parallel \nn \\
&\leq& \parallel \ \parallel B \parallel_2 \ \parallel u_a \parallel_{\infty} \ + \ \parallel R_1 \parallel_2 \ \parallel_1 \nn \\
&\leq& \left ( 2c\ N_2\ t^{-1/2} + r_1 \right ) h(t)
\eea

\noi so that
\beq
\label{2.28e}
N'_0 \leq 2c\ N_2\ T^{-1/2} + r_1 \ .
\eeq

\noi We next estimate $N'_1$. By Lemma 2.1
\bea
\label{2.29e}
\parallel v';L^{8/3}(J, L^4) \parallel &\leq& C \Big ( \parallel A_a v'; L^1(J,L^2)\parallel\ + \ \parallel Bv';L^{8/5} (J, L^{4/3})\parallel \nn \\
&&+ \parallel B u_a - R_1;L^1(J, L^2) \parallel\Big )\ .
\eea

\noi The last norm has already been estimated by (\ref{2.27e}), while
$$\parallel A_a v'; L^1(J, L^2)\parallel\ \leq \ \parallel\ \parallel A_a \parallel_{\infty} \ \parallel v'\parallel_2\ \parallel_1 \ \leq 4a \ N'_0 \ h(t)\ ,$$
$$\parallel B v'; L^{8/5}(J, L^{4/3})\parallel\ \leq \ \parallel\ \parallel B \parallel_{2} \ \parallel v'\parallel_4\ \parallel_{8/5} \ \leq C \ N_2\ N'_1 \ \overline{h}(t)\ h(t)$$

\noi by Lemma 2.2. Substituting those estimates into (\ref{2.29e}) yields 
$$N'_1 \leq C \left ( a\ N'_0 + c\ N_2\ T^{-1/2} + r_1 + N_2\ N'_1\ \overline{h}(T)\right )$$

\noi and therefore
\beq
\label{2.30e}
N'_1 \leq C_1 \left ( a\ N'_0 + c\ N_2\ T^{-1/2} + r_1 \right )
\eeq

\noi for $T$ sufficiently large satisfying a condition of the type $N_2 \overline{h}(T) \leq C$. We next estimate $N'_3$. From the time derivative of the equation for $v'$, we obtain
\bea
\label{2.31e}
&&\parallel \partial_t v'(t) \parallel_2^2 \ \leq \ \parallel \partial_t v'(t_0) \parallel_2^2\ + \ 2\parallel \ \parallel \partial_t v'\parallel_2 \Big ( \parallel \partial_t A_a \parallel_{\infty} \ \parallel v' \parallel_2\nn \\
&&+\ \parallel \partial_t B \parallel_2\ \parallel u_a \parallel_{\infty}\ + \ \parallel B \parallel_2\ \parallel \partial_t u_a \parallel_{\infty} \ + \ \parallel \partial_t R_1 \parallel_2 \Big ) \nn \\
&&+\ \parallel \partial_t v' \parallel_4\ \parallel \partial_t B \parallel_2\ \parallel v' \parallel_4\ \parallel_1
\eea

\noi We estimate the initial condition by 
\begin{eqnarray*}
\parallel \partial_t v'(t_0) \parallel_2 &\leq & \parallel B(t_0) \parallel_2\ \parallel u_a (t_0) \parallel_{\infty} \ + \ \parallel R_1(t_0) \parallel_2 \\
&\leq& \left ( c\ N_2\ t_0^{-3/2} + r_1 \right ) h(t_0)
\end{eqnarray*}

\noi where we have used the pointwise estimate
$$\parallel R_1(t) \parallel_2\ \leq\ \parallel \partial_t R_1;L^1([t, \infty ), L^2)\parallel \ \leq r_1\ h(t)\ .$$

\noi Using Lemma 2.2, we then obtain
$$N'^2_3 \leq \left ( c\ N_2\ T^{-3/2} + r_1 \right )^2 + C \left ( N'_3 \left ( a \ N'_0 + c\ N_2\ T^{-1/2} + r_1\right ) + N_2 \ N'_4\ N'_1\ \overline{h}(T)\right )$$

\noi and therefore
\beq
\label{2.32e}
N'_3 \leq C_3 \left ( a\ N'_0 + c\ N_2\ T^{-1/2} + r_1 + \left ( N_2\ N'_4\ N'_1\ \overline{h}(T)\right )^{1/2} \right ) \ .
\eeq

\noi We next estimate $N'_4$. By Lemma 2.1
\bea
\label{2.33e}
&&\parallel \partial_t v';L^{8/3}(J,L^4) \parallel \ \leq \ C \Big ( \parallel \ \parallel A_a \parallel_{\infty} \ \parallel \partial_t v'\parallel_2 \ + \  \parallel \partial_t A_a \parallel_{\infty} \ \parallel v' \parallel_2\nn \\
&&+\ \parallel \partial_t B \parallel_2\ \parallel u_a \parallel_{\infty}\ + \ \parallel B \parallel_2\ \parallel \partial_t u_a \parallel_{\infty} \ + \ \parallel \partial_t R_1 \parallel_2 \ \parallel_1\nn \\
&&+\ \parallel\ \parallel \partial_t B \parallel_2\ \parallel v' \parallel_4\ + \ \parallel B \parallel_2\ \parallel \partial_t v' \parallel_4\ \parallel_{8/5}\Big )
\eea

\noi and therefore by Lemma 2.2
\beq
\label{2.34e}
N'_4 \leq C_4 \left ( a\left ( N'_3 + N'_0\right ) + c\ N_2\ T^{-1/2} + r_1 + N_2\left (  N'_1 +  N'_4\right )\ \overline{h}(T)\right )\ .
\eeq

\noi We next estimate $N'_5$. From (\ref{1.5e}) we obtain directly for $2 \leq r \leq 4$
\bea
\label{2.35e}
\parallel \Delta v' \parallel_r &\leq& 2 \Big ( \parallel \partial_t v' \parallel_r\ + \ \parallel A_a \parallel_{\infty}\ \parallel v' \parallel_r \ + \ \parallel B \parallel_r \ \parallel u_a \parallel_{\infty}\nn \\
&&+ \ \parallel R_1 \parallel_r \ + \ C \parallel B \parallel_r\ \parallel v' \parallel_r^{1-3/2r}\ \parallel \Delta v' \parallel_r^{3/2r}\Big )
\eea

\noi by a Sobolev inequality, and therefore for $r=2$
\beq
\label{2.36e}
N'_5 \leq 4 \left ( N'_3 + a\ N'_0 \ T^{-1} + c\ N_2\ T^{-3/2} + r_1 + C\ N'_0 \left ( N_2\ h(T)\right )^4\right ) \ .
\eeq

\noi Similarly, taking the norm in $L^{8/3}(J)$ of (\ref{2.35e}) with $r = 4$, we obtain 
\beq
\label{2.37e}
N'_6 \leq 4 \left ( N'_4 + a\ N'_1 \ T^{-1} + C\ c\ N_2\ T^{-9/8} + r_1 \ T^{-\eta} + C\ N'_1 \left ( N_2\ h(T)\right )^{8/5}\right ) 
\eeq

\noi where we have used the estimate
$$\parallel B \parallel_4\ \leq\ C \parallel B \parallel_2^{1/4} \ \parallel \nabla B \parallel_2^{3/4}\ \leq \ C\ N_2\ h(t)\ .$$

\noi We next turn to the estimate of $N'_2$. The equation for $B'$ takes the form $\sq B' = \Delta F$ where
$$F = |v|^2 + 2 {\rm Re} \ \overline{u}_a v - \Delta^{-1} R_2$$

\noi and by standard energy estimates
\beq
\label{2.38e}
\left \{ \begin{array}{l}\parallel \nabla B'(t) \parallel_2\ \vee\ \parallel \partial_t B'(t) \parallel_2 \ \leq \ \parallel \Delta F; L^1(J, L^2)\parallel\\ \\
\parallel B'(t) \parallel_2\ \vee\ \parallel \omega^{-1} \partial_t B'(t) \parallel_2 \ \leq \ \parallel \nabla F; L^1(J, L^2)\parallel \ .
\end{array}\right .
\eeq

\noi We estimate
\bea
\label{2.39e}
\parallel \nabla F;L^1(J,L^2) \parallel &\leq& \parallel 2\ \parallel v \parallel_4\ \parallel \nabla v \parallel_4 \ + \ \parallel v \parallel_2\ \parallel \nabla u_a\parallel_{\infty}\nn \\
&+& \parallel \nabla v \parallel_2\ \parallel u_a \parallel_{\infty} \ + \ \parallel \omega^{-1} R_2 \parallel_2\ \parallel_1  
\eea
\bea
\label{2.40e}
&&\parallel \Delta F;L^1(J,L^2) \parallel \ \leq \  \parallel 2\left (  \parallel v \parallel_4\ \parallel \Delta v \parallel_4 \ + \ \parallel \nabla v \parallel_4^2\right )\ + \  \parallel v \parallel_2\ \parallel \Delta  u_a\parallel_{\infty}\nn \\
&&+\ 2 \parallel \nabla v \parallel_2\ \parallel \nabla u_a \parallel_{\infty} \ + \ \parallel \Delta v \parallel_2\ \parallel u_a \parallel_{\infty} \ + \ \parallel R_2 \parallel_2\ \parallel_1\ .
\eea

\noi Using Lemma 2.2 and the definitions, we obtain from (\ref{2.38e})-(\ref{2.40e})
\beq
\label{2.41e}
N'_2 \leq C_2 \left ( c (N_0 + N_5)T^{-1/2} + N_1(N_1 + N_6) \overline{h}(T) + r_2 \right ) \ .
\eeq

It follows immediately from (\ref{2.28e}) (\ref{2.30e}) (\ref{2.32e})
(\ref{2.34e}) (\ref{2.36e}) (\ref{2.37e}) (\ref{2.41e}) that $(v', B')$
is bounded in $X(I)$ uniformly in $t_0$ for $T$ sufficiently large,
more precisely for $N_2 \overline{h}(T) \leq C$ for a suitable absolute
constant $C$. \par

From now on the proof is very similar to that of Proposition 2.2 in
\cite{11r}. We next take the limit $t_0 \to \infty$ of $(v'_{t_0},
B'_{t_0})$, restoring the subscript $t_0$ for that part of the
argument. Let $T < t_0 < t_1 < \infty$ and let $(v'_{t_0}, B'_{t_0})$
and $(v'_{t_1}, B'_{t_1})$ be the corresponding solutions of
(\ref{1.5e}). From the $L^2$ norm conservation of the difference
$v'_{t_0}- v'_{t_1}$ and from (\ref{2.28e}), it follows that for all $t
\in [T, t_0]$ 
\beq
\label{2.42e}
\parallel v'_{t_0}(t)- v'_{t_1}(t)\parallel_2\ = \ \parallel v'_{t_1}(t_0)\parallel_2\ \leq \ K_0\ h(t_0)
\eeq 

\noi where $K_0$ is the RHS of (\ref{2.28e}), while from (\ref{1.5e})
(\ref{2.38e})-(\ref{2.41e}) and the initial
conditions, it follows that
\beq
\label{2.43e}
\parallel B'_{t_0} - B'_{t_1};L^{\infty}([T, t_0], H^1)\parallel\ \vee \ \parallel \partial_t (B'_{t_0} - B'_{t_1});L^{\infty}([T, t_0],L^2)\parallel \ 
\leq \  K_2 \  h(t_0)
\eeq

\noi where $K_2$ is the RHS of (\ref{2.41e}).\par

It follows from (\ref{2.42e}) (\ref{2.43e}) that there exists $(v', B')
\in L_{loc}^{\infty}([T, \infty ), L^2 \oplus H^1)$ with $\partial_tB' \in L_{loc}^{\infty}([T, \infty ), L^2)$ such that
$(v'_{t_0}, B'_{t_0})$ converges to $(v', B')$ in that space when $t_0
\to \infty$. From the uniformity in $t_0$ of the estimates
(\ref{2.28e}) (\ref{2.41e}), it follows that $(v', B')$
satisfies the same estimates in $[T, \infty )$, namely that
(\ref{2.28e}) (\ref{2.41e}) hold with $N'_i$ defined by
(\ref{2.20e}) (\ref{2.22e}) with $I = [T, \infty )$.
Furthermore it follows by a standard compactness  argument that $(v',
B') \in X([T, \infty ))$ and that $v'$ satisfies the remaining
estimates, namely (\ref{2.30e}) (\ref{2.32e}) (\ref{2.34e})
(\ref{2.36e}) (\ref{2.37e}) with the remaining $N'_i$ again defined by (\ref{2.21e})
(\ref{2.23e})-(\ref{2.26e}) with $I = [T, \infty )$.
Clearly $(v', B')$ satisfies the system (\ref{1.5e}).

From now on, $I$ denotes the interval $[T, \infty )$. The previous
construction defines a map $\phi : (v, B) \to (v', B')$ from $X(I)$ to
itself. The next step consists in proving that the map $\phi$ is a
contraction on a suitable closed bounded set ${\cal R}$ of $X(I)$. We
define ${\cal R}$ by the conditions (\ref{2.13e})-(\ref{2.19e}) for
some constants $N_i$ and for all $t \in I$. We first show that for a
suitable choice of $N_i$ and for sufficiently large $T$, the map $\phi$
maps ${\cal R}$ into ${\cal R}$. By (\ref{2.28e}) (\ref{2.30e})
(\ref{2.32e}) (\ref{2.34e}) (\ref{2.36e}) (\ref{2.37e}) (\ref{2.41e}), it suffices for that purpose that

\beq \label{2.44e} \left \{ \begin{array}{l} 
(N'_0 \leq ) \ r_1 + 2c  \ N_2 T^{-1/2} \leq N_0 \\ \\ (N'_1 \leq ) \ C_1 \Big ( r_1 + a N'_0 + c\ N_2  T^{-1/2}\Big )  \leq N_1 \\ \\ 
  (N'_2 \leq ) \ C_2 \Big ( r_2 + c(N_0 + N_5) T^{-1/2} + N_1(N_1+N_6) \overline{h}(T)\Big ) \leq N_2 \\ \\ 
(N'_3 \leq ) \ C_3 \Big ( r_1 + a N'_0 + cN_2 T^{-1/2} + (N_2 N'_4 N'_1 \overline{h}(T))^{1/2} \Big )\leq N_3 \\ \\ 
(N'_4 \leq ) \ C_4 \Big ( r_1 + a(N'_3 + N'_0) + c N_2 T^{-1/2} + N_2 (N'_1 + N'_4) \overline{h}(T) \Big ) \leq N_4\\ \\ 
(N'_5 \leq ) \ 4 \Big ( r_1 + N'_3 + a N'_0 T^{-1} + c N_2 T^{-3/2} + C N'_0 (N_2h(T))^4 \Big ) \leq N_5\\ \\ 
(N'_6
\leq ) \ 4 \Big ( r_1 + N'_4 + a N'_1 T^{-1} + C c N_2 T^{-9/8} + C N'_1 (N_2 h(T))^{8/5}\Big ) \leq N_6\ . \end{array}
\right . \eeq

We ensure those conditions as follows. We ensure the first two conditions by taking
\beq
\label{2.45e}
\left \{ \begin{array}{l} N_0 = r_1 + 1 \\ \\ N_1 = C_1 \left ( r_1 + aN_0 + 1 \right ) 
\end{array}
\right . \eeq

\noi and by taking $T$ sufficiently large for the $o(1)$ terms in those
conditions not to exceed 1. It is then easy to see that the conditions
on $N_3$, $N_4$ are satisfied by taking
\beq
\label{2.46e}
\left \{ \begin{array}{l} N_3 = C_3 \left ( r_1 + aN_0 + 1\right )  \\ \\ N_4 = C_4 \left ( r_1 + a(N_3 + N_0) + 1 \right ) 
\end{array}
\right . \eeq

\noi and by taking $T$ sufficiently large for the $o(1)$ terms in those
conditions with the $N'_i$ replaced by $N_i$ not to exceed 1. We
finally take
\beq
\label{2.47e}
\left \{ \begin{array}{l} N_5 = 4 \left ( r_1 + N_3 + 1\right )  \\ \\ N_6 = 4 \left ( r_1 + N_4 + 1 \right )\\ \\ N_2 = C_2\left ( r_2 + 1\right ) 
\end{array}
\right . \eeq

\noi and we take in addition $T$ sufficiently large to ensure that the
$o(1)$ terms in the corresponding conditions do not exceed 1. This
completes the proof of the stability of ${\cal R}$.

 We next show that the map $\phi$ is a contraction on ${\cal R}$.
Let $(v_i, B_i) \in {\cal R}$, $i = 1,2$, and let $(v'_i, B'_i) = \phi
((v_i, B_i))$. For any pair of functions $(f_1, f_2)$ we define $f_{\pm}
= (1/2)(f_1 \pm f_2)$ so that $(fg)_{\pm} = f_+ g_{\pm} + f_-
g_{\mp}$. In particular $u_+ = u_a + v_+$, $u_- = v_-$, $A_+ = A_a +
B_+$ and $A_- = B_-$. Corresponding to (\ref{1.5e}), $(v'_-, B'_-)$
satisfies the system
\beq
\label{2.48e}
\left \{ \begin{array}{l} i\partial_t v'_- = - (1/2) \Delta v'_- + A_+ v'_- + B_- u_a + B_- v'_+\\ \\ \sq B'_- =  2\Delta\  {\rm Re}\ \left ( \overline{u}_a + \overline{v}_+ \right ) v_-\ .\end{array} \right .
\eeq

\noi Since ${\cal R}$ is convex and stable under $\phi$, $(v_+, B_+)$
and $(v'_+, B'_+)$ belong to ${\cal R}$, namely satisfy
(\ref{2.13e})-(\ref{2.19e}). Let $N_{i^-}$ and $N'_{i^-}$ be the seminorms of $(v_-, B_-)$ and $(v'_-, B'_-)$ corresponding to
(\ref{2.20e})-(\ref{2.26e}), namely the constants obtained by replacing
$(v', B', N'_i)$ by $(v_-, B_-, N_{i^-})$ and $(v'_-, B'_-, N'_{i^-})$
in (\ref{2.20e})-(\ref{2.26e}). We have to estimate the $N'_{i^-}$ in
terms of the $N_{i^-}$. The estimates are essentially the same as those
of $N'_i$ in terms of $N_i$ with minor differences~: the contribution
of the remainders disappear, the linear terms are the same, and the
quadratic terms are in general obtained by polarization. The only
exceptions to that rule are the $B_- v'_+$ term in the estimate of
$N'_{0^-}$ and the $B_-\partial_t v'_+$ term in the estimate of
$N'_{3^-}$ because the corresponding terms in the estimate of one single function disappear for algebraic reasons. Thus we estimate
\bea
\label{2.49e}
&&\parallel v'_-(t) \parallel_2^2\ \leq \ 2 \parallel \ <v'_-, B_-(u_a + v'_+) >\   \parallel_1\nn \\
&&\leq \ 2 \parallel \ \parallel v'_- \parallel_2\  \parallel B_- \parallel_2\  \parallel u_a \parallel_{\infty} \ + \  \parallel v'_- \parallel_4\  \parallel B_- \parallel_2\  \parallel v'_+ \parallel_4 \  \parallel_1
\eea

\noi and therefore by Lemma 2.2
$$N'^2_{0^-} \leq 2c N'_{0^-} \ N_{2^-} \ T^{-1/2} + C\ N'_{1^-}\ N_{2^-}\ N_1 \ \overline{h}(T)$$

\noi so that 
\beq
\label{2.50e}
N'_{0^-} \leq 2c \ N_{2^-} \ T^{-1/2} + C \left ( N'_{1^-}\ N_{2^-}\ N_1\ \overline{h}(T)\right )^{1/2}\ .
\eeq

\noi Similarly
\bea
\label{2.51e}
&&\parallel \partial_t v'_- (t)\parallel_2^2\ \leq\ 2 \parallel \ \parallel \partial_t v'_-) \parallel_2 \Big (  \parallel \partial_tA_a \parallel_{\infty} \  \parallel v'_- \parallel_2 \ + \  \parallel \partial_t B_-\parallel_2\  \parallel u_a \parallel_{\infty}\nn\\
&&+ \ \parallel B_-\parallel_2  \ \parallel \partial_t u_a \parallel_{\infty}\Big ) \ + \  \parallel \partial_t v'_- \parallel_4 \Big (  \parallel \partial_t B_+ \parallel_2\  \parallel v'_- \parallel_4\ + \  \parallel \partial_t B_- \parallel_2\  \parallel v'_+ \parallel_4\nn\\
&&+ \ \parallel B_- \parallel_2\  \parallel \partial_t v'_+ \parallel_4 \Big ) \parallel_1
\eea

\noi and therefore by Lemma 2.2
\beq
\label{2.52e}
N'^2_{3^-} \leq C \left ( N'_{3^-} \left ( a N'_{0^-} + c N_{2^-} T^{-1/2}\right ) + N'_{4^-} \left ( N_2 N'_{1^-} + N_{2^-} (N_1 + N_4)\right ) \overline{h}(T)\right )
\eeq
\beq
\label{2.53e}
N'_{3^-} \leq C_3 \left ( a N'_{0^-}  + c N_{2^-} T^{-1/2} + \left ( N'_{4^-} \left ( N_2 N'_{1^-} + N_{2^-} (N_1 + N_4)\right ) \overline{h}(T)\right )^{1/2} \right ) \ .
\eeq

\noi The estimates of the other $N'_{i^-}$ follow the general rule and are thus given by
\beq
\label{2.54e}
N'_{1^-} \leq C_1 \left ( a N'_{0^-} + c N_{2^-} T^{-1/2} + N_{2^-} N_1 \overline{h}(T)\right )
\eeq
\beq
\label{2.55e}
N'_{4^-} \leq C_4 \left ( a \left ( N'_{3^-} + N'_{0^-}\right )   + c N_{2^-} T^{-1/2} + \left (  N_2 N'_{1^-} + N_{2^-} (N_1 + N_4)\right ) \overline{h}(T)\right ) 
\eeq
\bea
\label{2.56e}
&&N'_{5^-} \leq 4 \Big ( N'_{3^-} + a N'_{0^-} T^{-1} + c N_{2^-} T^{-3/2} + C N'_{0^-}\left ( N_2 h(T)\right )^4 \nn \\
&&+ C N_{2^-} N_0^{1/4} N_5^{3/4} h(t) \Big ) 
\eea
\bea
\label{2.57e}
&&N'_{6^-} \leq 4 \Big ( N'_{4^-} + a N'_{1^-} T^{-1} + C c N_{2^-} T^{-9/8} + C N'_{1^-}\left ( N_2 h(T)\right )^{8/5}\nn \\
&& + \ C N_{2^-} N_1^{5/8} N_6^{3/8} h(t) \Big ) 
\eea

\noi where in the last term we have estimated
$$\parallel B_- v'_+ \parallel_4\ \leq \ C \parallel B_- \parallel_4 \ \parallel v'_+ \parallel_4^{5/8}\ \parallel \Delta v'_+ \parallel_4^{3/8} \ .$$

\noi Finally,
\beq
\label{2.58e}
N'_{2^-} \leq C_2 \left ( c \left ( N_{0^-} + N_{5^-}\right ) T^{-1/2} + \left ( N_{1^-} (N_1 + N_6) + N_{6^-} N_1 \right ) \overline{h} (T)\right ) \ .
\eeq

We have kept the same constants $C_i$ in (\ref{2.54e}) (\ref{2.53e})
(\ref{2.55e}) (\ref{2.58e}) as in (\ref{2.30e}) (\ref{2.32e})
(\ref{2.34e}) (\ref{2.41e}). In fact those constants are determined by
the linear terms in the estimates, which are the same in both cases.
There may occur additional different constants coming from the
quadratic terms. They have been omitted in
(\ref{2.53e})-(\ref{2.58e}).\par

From the fact that most of the terms in the RHS of (\ref{2.50e}) and of
(\ref{2.53e})-(\ref{2.58e}) are $o(1)$ when $T \to \infty$ and that
this system of inequalities is strictly triangular in the $o(1)$
terms, it follows easily as in \cite{11r} \cite{12r} that the map
$\phi$ is a contraction in the set of semi norms $N_i$, $0 \leq i \leq
6$, for $T$ sufficiently large. It follows therefrom that the system
(\ref{1.3e}) has a unique solution in ${\cal R}$. Uniqueness in
$X(\cdot )$ follows from the same estimates. \par

The last statement of the Proposition follows from the estimates of
$B'$ leading to (\ref{2.41e}) (see especially
(\ref{2.38e})-(\ref{2.40e})) by using the stronger estimate
(\ref{2.11e}) of $R_2$.\par\nobreak \hfill $\sq$ \par

We now turn to the second step of the method, namely to the choice of
$(u_a, A_a)$ and the derivation of the conditions
(\ref{2.6e})-(\ref{2.11e}). We shall need the standard factorisation of
the free Schr\"odinger group
\beq
\label{2.59e}
U(t) \equiv \exp (i(t/2) \Delta ) = M\ D\ F\ M
\eeq
\noi where
\beq
\label{2.60e}
M \equiv M(t) = \exp (ix^2/2t)
\eeq
\beq
\label{2.61e}
D(t) = (it)^{-n/2} D_0(t) \quad , \quad \left ( D_0(t) f \right ) (x) = f(x/t) \ .
\eeq

\noi Using that decomposition, one can easily derive the following lemma, which we state for $n = 2, 3$. \\

\noi {\bf Lemma 2.3.} {\it  Let $n = 2$ or $3$. Let $u_+ \in
H^{0,2}(\subset L^1)$ and let $u_0 = U(t) u_+$. Then the following estimates
hold~:}
\beq
\label{2.62e}
\parallel \nabla |u_0|^2\parallel_2 \ \leq 2 (2 \pi t)^{-n/2} \ t^{-1} \parallel u_+\parallel_1 \ \parallel xu_+\parallel_2 \ ,
\eeq
\beq
\label{2.63e}
\parallel \Delta |u_0|^2\parallel_2 \ \leq 4 (2 \pi t)^{-n/2} \ t^{-2} \parallel u_+\parallel_1 \ \parallel x^2u_+\parallel_2 \ .
\eeq
\vskip 5 truemm

\noi {\bf Proof.} From the representation (\ref{2.59e}), we obtain
\begin{eqnarray*}
\parallel \nabla |u_0|^2 \parallel_2 &\leq& 2 t^{-n-1} \parallel D_0(t) (\overline{FMu_+} FMxu_+) \parallel_2 \\
&\leq& 2 t^{-n/2-1} \parallel FMu_+\parallel_{\infty} \ \parallel FMxu_+\parallel_2\\
&\leq& 2 (2 \pi t)^{-n/2} t^{-1} \parallel u_+\parallel_1\ \parallel xu_+ \parallel_2 \end{eqnarray*}

\noi by the Hausdorf-Young inequality. \par

Similarly
$$\parallel \Delta |u_0|^2 \parallel_2 \ \leq \ 2 t^{-n-2} \left ( \parallel D_0(t) (\overline{FMu_+} FMx^2u_+) \parallel_2 \ + \ \parallel D_0(t) |FMxu_+|^2 \parallel_2 \right )$$
$$\leq \ 2 t^{-n/2-2} \left ( \parallel FMu_+ \parallel_{\infty} \ \parallel FM x^2 u_+ \parallel_2 \ + \ \parallel FMx u_+ \parallel_4^2 \right )$$
$$\leq \ 4 (2 \pi t)^{-n/2}t^{-2} \parallel u_+ \parallel_1\ \parallel x^2u_+ \parallel_2$$

\noi by the Hausdorf-Young and H\"older inequalities.\par \nobreak \hfill $\sq$ \par

We shall also need some estimates of solutions of the free wave
equation, which we collect in the following lemma. A proof can be found
in \cite{24r}.\\

\noi {\bf Lemma 2.4.} {\it Let $A_0$ be defined by (\ref{1.8e}). Let $k \geq 0$ be an integer. Let $A_+$ and $\dot{A}_+$ satisfy the conditions
\beq
\label{2.64e}
A_+, \omega^{-1}\dot{A}_+ \in H^k \qquad , \quad \nabla^2 A_+, \nabla \dot{A}_+ \in W_1^k \ .
\eeq
\noi Then $A_0$ satisfies estimates
\beq
\label{2.65e}
\left \{ \begin{array}{l}  \parallel A_0(t);W_r^k \parallel \ \leq a\ t^{-1+2/r} \ ,\\ \\ \parallel \partial_t A_0(t);W_r^{k-1}\parallel\ \leq a\ t^{-1+2/r} \quad \hbox{\it for $k \geq 1$}\ .\end{array} \right .
\eeq

\noi for $2 \leq r \leq \infty$ and for all $t\in {I\hskip-1truemm R}$, where $a$
depends on $A_+$, $\dot{A}_+$ through the norms associated with
(\ref{2.64e}).}\\

We are now in a position to derive the final result with simple
asymptotics (\ref{1.7e}) (\ref{1.8e}), namely Proposition 1.1.\\

\noi {\bf Proof of Proposition 1.1.}\par

The result will follow from Proposition 2.1 once we have proved that
$(u_0, A_0)$ satisfies the assumptions of that proposition for $(u_a,
A_a)$. From the standard $L^1-L^{\infty}$ estimates of $U(t)$, we obtain 
\beq
\label{2.66e}
\parallel u_0 (t) \parallel_{\infty} \ \leq \ (2 \pi t)^{-3/2} \parallel u_+ \parallel_1\ ,
\eeq
\beq
\label{2.67e}
2 \parallel \partial_t u_0 \parallel_{\infty}\ = \ \parallel \Delta u_0 \parallel_{\infty} \ \leq \ (2 \pi t)^{-3/2} \parallel \Delta u_+ \parallel_1 \ ,
\eeq

\noi which proves (\ref{2.6e}). The assumption (\ref{2.7e}) on $A_0$
follows from Lemma 2.4. We next consider the remainders $R_1 = - A_0
u_0$ and $R_2 = - \Delta |u_0|^2$. We estimate
$$\parallel R_1 \parallel_2 \ \leq \ \parallel A_0 \parallel_2\ \parallel u_0 \parallel_{\infty} \ \leq \ C\ t^{-3/2}$$
$$\parallel \partial_t R_1 \parallel_2 \ \leq \ \parallel A_0 \parallel_2\ \parallel \partial_t u_0 \parallel_{\infty} \ + \ \parallel \partial_t A_0 \parallel_2 \ \parallel u_0 \parallel_{\infty} \ \leq C \ t^{-3/2}$$

\noi by (\ref{2.66e}) (\ref{2.67e}) and Lemma 2.4. This proves (\ref{2.8e}) with $h(t) = t^{-1/2}$. On the other hand 
$$\parallel R_1 \parallel_4 \ \leq \ \parallel A_0 \parallel_4\ \parallel u_0 \parallel_{\infty} \ \leq C\ t^{-2}$$

\noi by (\ref{2.66e}) and Lemma 2.4, which yields (\ref{2.9e}) with $\eta = 9/8$. Finally
$$\parallel R_2 \parallel_2 \ = \ \parallel \Delta |u_0|^2 \parallel_2\ \leq \ 2 \left ( \parallel u_0 \parallel_{\infty} \ \parallel \Delta u_0 \parallel_2 \ + \ \parallel \nabla u_0 \parallel_4^2 \right ) \ \leq C\ t^{-3/2}\ ,$$
$$\parallel \omega^{-1} R_2 \parallel_2 \ = \ \parallel \nabla |u_0|^2\parallel_2\ \leq\ 2 \parallel u_0 \parallel_{\infty} \ \parallel \nabla u_0 \parallel_2\ \leq \ C\ t^{-3/2} \ ,$$

\noi which proves (\ref{2.10e}).\par

The last statement of the proposition follows from the corresponding statement of Proposition 2.1 and from Lemma 2.3 with $n = 3$, which yields actually
$$\parallel \omega^{-1} R_2; H^1 \parallel\ \leq C\ t^{-5/2}$$

\noi and therefore
$$\parallel \omega^{-1}R_2; L^1([t, \infty ), H^1) \parallel \ \leq C \ t^{-1}\ h(t)$$

\noi which is stronger than (\ref{2.11e}) by a factor $t^{-1/2}$.\par \nobreak \hfill $\sq$ \par

We next turn to the case where one uses the more accurate asymptotic
form proposed in \cite{23r}, thereby obtaining a stronger asymptotic
convergence in time of the solution on a smaller subspace of asymptotic
states. Thus we choose
\beq
\label{2.68e}
(u_a, A_a) = \left ( (1 + f) u_0, A_0 \right )
\eeq

\noi where $(u_0, A_0)$ are defined by (\ref{1.7e}) (\ref{1.8e}) and
\beq
\label{2.69e}
f = 2 \Delta^{-1} A_0 \ .
\eeq 

\noi Using the operators
\beq
\label{2.70e}
J = x + it \nabla \quad , \quad P = t \partial_t + x \cdot \nabla \ ,
\eeq

\noi we can rewrite the remainders $R_1$ and $R_2$ as
\bea
\label{2.71e}
R_1 &=& \left ( i \partial_t + (1/2) \Delta - A_0\right ) (1 + f) u_0\nn \\
&=& - f\ A_0\ u_0 - it^{-1} (\nabla f)\cdot Ju_0 + i t^{-1} (Pf) u_0
\eea
\beq
\label{2.72e}
R_2 = - \Delta (1 + f)^2|u_0|^2\ .
\eeq

\noi We first reduce the estimates required for $R_1$ and $R_2$ to general estimates of $u_+$, $A_0$ and $f$. We first estimate $R_1$.\\

\noi {\bf Lemma 2.5.} {\it Let $u_+ \in W_1^2$, $xu_+ \in W_1^2$, and let $A_0$ and $f$ satisfy 
\beq
\label{2.73e}
\parallel \partial_t^j \nabla^k A_0\parallel_{\infty} \ \leq a\ t^{-1} 
\eeq
\beq
\label{2.74e}
\parallel \partial_t^j \nabla^k f;H^1\parallel\ \vee \  \parallel \partial_t^j \nabla^k Pf\parallel_2\ \leq C 
\eeq

\noi for $0 \leq j + k \leq 1$ and for all $t \geq 1$. Then the following estimates hold~:
\beq
\label{2.75e}
\parallel \partial_t^j \nabla^k R_1\parallel_2 \ \leq C\ t^{-5/2} 
\eeq

\noi for some constant $C$, for $0 \leq j + k \leq 1$ and for all $t \geq 1$.}\\

\noi {\bf Proof.} By the $L^1-L^{\infty}$ estimate of $U(t)$ and the commutation rule $JU(t) = U(t)x$, we obtain 
$$\parallel \partial_t^j \nabla^k u_0 \parallel_{\infty} \ \vee \ \parallel \partial_t^j \nabla^k J u_0 \parallel_{\infty} \ \leq C\ t^{-3/2}$$

\noi for $0 \leq j + k \leq 1$. We then estimate
\begin{eqnarray*}
\parallel R_1 \parallel_2 &\leq& \parallel f \parallel_2 \ \parallel A_0 \parallel_{\infty} \ \parallel u_0 \parallel_{\infty} \ + \ t^{-1} \parallel \nabla f \parallel_2\ \parallel J u_0 \parallel_{\infty}\\
&&+ \ t^{-1}  \parallel Pf \parallel_2\ \parallel u_0 \parallel_{\infty} \ \leq C\ t^{-5/2}
\end{eqnarray*}

\noi which proves (\ref{2.75e}) for $j = k = 0$. The other cases are
obtained similarly by distributing $\partial_t$ or $\nabla$ among the
various factors. \par\nobreak \hfill $\sq$\par

We next estimate $R_2$.\\

\noi {\bf Lemma 2.6.} {\it Let $u_+ \in W_1^2 \cap H^{0,2}$ and let $f$ satisfy 
\beq
\label{2.76e}
\parallel \nabla f(t)\parallel_2\ \vee \  \parallel \Delta f(t)\parallel_2\ \vee \ \parallel f(t) \parallel_{\infty}\ \leq C 
\eeq

\noi for all $t \geq 1$. Then the following estimates hold~: 
\beq
\label{2.77e}
\parallel \omega^{-1}R_2 \parallel_2 \ \leq C\ t^{-5/2} \ ,
\eeq
\beq
\label{2.78e}
\parallel R_2 \parallel_2 \ \leq C\ t^{-3} 
\eeq

\noi for some constant $C$ and for all $t \geq 1$.}\\

\noi {\bf Proof.} For $u_+ \in W_1^2$, we know that $\parallel \nabla^j
u_0 \parallel_{\infty} \leq C t^{-3/2}$ and therefore $\parallel
\nabla^j|u_0|^2 \parallel_{\infty} \leq C t^{-3}$ for $j = 0, 1, 2$.
For $u_+ \in H^{0,2}$, we know that $\parallel \nabla |u_0|^2
\parallel_2 \leq C t^{-5/2}$ and $\parallel \Delta |u_0|^2 \parallel_2
\leq C t^{-7/2}$ by Lemma 2.3. We then estimate $R_2$ as follows
$$\parallel \omega^{-1} R_2\parallel_2\ = \ \parallel \nabla (1 + f)^2|u_0|^2 \parallel_2$$
$$\leq \ \left ( 1 \ + \ \parallel f\parallel_{\infty}\right )^2 \ \parallel \nabla |u_0 |^2 \parallel_2\ + \ 2 \left ( 1 \ + \ \parallel f \parallel_{\infty}\right ) \ \parallel \nabla f \parallel_2\ \parallel u_0 \parallel_{\infty}^2$$
$$\leq C\ t^{-5/2}\ , $$

$$\parallel R_2\parallel_2\ = \ \parallel \Delta (1 + f)^2|u_0|^2 \parallel_2$$
$$\leq \ \left ( 1 \ + \ \parallel f\parallel_{\infty}\right )^2 \ \parallel \Delta |u_0 |^2 \parallel_2\ + \ 4 \left ( 1 \ + \ \parallel f \parallel_{\infty}\right ) \ \parallel \nabla f \parallel_2\ \parallel \nabla |u_0|^2 \parallel_{\infty}$$
$$+\ 2 \left ( \left ( 1\ + \ \parallel f \parallel_{\infty}\right ) \parallel \Delta f \parallel_2\ + \ \parallel \nabla f \parallel_4^2 \right ) \parallel u_0 \parallel_{\infty}^2\ \leq C\ t^{-3}\ . $$
\hfill $\sq$\par

\noi {\bf Remark 2.1.} In practice the bound on $\parallel f \parallel_{\infty}$ in (\ref{2.76e}) will follow from the Sobolev inequality
\beq
\label{2.79e}
\parallel f \parallel_{\infty}^2\ \leq \ C \parallel \nabla f \parallel_2 \ \parallel \Delta f \parallel_2
\eeq

\noi for $f$ tending to zero at infinity in some weak sense.\\

We are now in a position to derive the final result with improved
asymptotics (\ref{2.68e}) (\ref{2.69e}), namely Proposition 1.2. \\

\noi {\bf Proof of Proposition 1.2}\par

\noi \underbar{Part 1}. The result follows from Proposition 2.1 and
from the fact that $(u_a, A_0)$ satisfies the assumptions of that
proposition for $(u_a, A_a)$. The condition (\ref{2.6e}) for $u_a$
follows from the same condition for $u_0$, which follows from
(\ref{2.66e}) (\ref{2.67e}), and from $L^{\infty}$ estimates for $f$. For $\partial = 1$,
$\partial_t$, $\nabla$, we estimate
\beq
\label{2.80e}
\parallel \partial f \parallel_{\infty}^2\ \leq \ C \parallel \nabla \partial f \parallel_2^{1/2}  \ \parallel \Delta \partial f \parallel_2^{1/2}\ \leq C
\eeq

\noi by a Sobolev inequality, the definition (\ref{2.69e}) of $f$ and Lemma 2.4, while
\beq
\label{2.81e}
\parallel \Delta f \parallel_{\infty}\ = \ C \parallel A_0 \parallel_{\infty} \ \leq\ C\ t^{-1} 
\eeq

\noi as a special case of (\ref{2.7e}), which also follows from
(\ref{1.10e}) and from Lemma 2.4 as before. The conditions (\ref{2.8e})
(\ref{2.9e}) with $h(t) = t^{-3/2}$ follow from Lemma 2.5, especially
(\ref{2.75e}), under the assumptions made on $u_+$ and the conditions
(\ref{2.73e}) (\ref{2.74e}). The latter follow from (\ref{1.10e}), from
Lemma 2.4, from the definition (\ref{2.69e}) of $f$ and from the fact
that $Pf$ is a solution of the free wave equation with initial data $(2x
\cdot \nabla \Delta^{-1} A_+, 2(1 + x \cdot \nabla ) \Delta^{-1}
\dot{A}_+)$.  Finally the condition (\ref{2.10e}) follows from Lemma
2.6, from (\ref{2.69e}), from (\ref{1.10e}) and from Lemma 2.4.\\

\noi \underbar{Part 2}. The result follows from the fact that $(fu_0,0)
\in X(I)$, namely that $fu_0$ satisfies the conditions on $v$ that
appear in the definition of $X(I)$, as we now show. For $u_+ \in
W_1^2$, we estimate
\bea
\label{2.82e}
\parallel \partial^{\alpha} (fu_0)\parallel_r &\leq & \sum_{\beta \leq \alpha} \parallel \partial^{\beta} f \parallel_r \ \parallel \partial^{\alpha - \beta} u_0 \parallel_{\infty}\nn \\
&\leq& C\ t^{-3/2} \sum_{\beta \leq \alpha} \parallel \partial^{\beta} f \parallel_r
\eea

\noi for $2 \leq r \leq \infty$ and $\alpha$ a multiindex with $|\alpha | \leq 2$, and
\bea
\label{2.83e}
\parallel \partial_t (fu_0)\parallel_r &\leq & \parallel f\parallel_r \ \parallel\partial_t u_0 \parallel_{\infty}\ + \  \parallel \partial_t f\parallel_r\ \parallel u_0 \parallel_{\infty}\nn \\
&\leq& C\ t^{-3/2} \left ( \parallel f\parallel_r\ + \ \parallel \partial_t f \parallel_r\right ) \ .
\eea

\noi For $r = 2$, it follows from (\ref{2.82e}) (\ref{2.83e}), from
(\ref{1.10e}), from Lemma 2.4 and from the definition (\ref{2.69e}) of
$f$ that 
\beq
\label{2.84e}
\parallel \partial_t (fu_0)  \parallel_{2}\ \vee \ \parallel fu_0;H^2\parallel \  \leq C \ t^{-3/2}\ . 
\eeq

\noi For $r=4$, it follows from the standard $L^p-L^q$ estimates for
the wave equation \cite{24r} and from the definition of $f$ that
\beq
\label{2.85e}
\parallel f\parallel_4\ \leq\ C\ t^{-1/2} \left ( \parallel \omega^{-1} A_+ \parallel_{4/3}\ + \ \parallel \omega^{-2} \dot{A}_+\parallel_{4/3} \right )
\eeq
\beq
\label{2.86e}
\parallel \nabla f\parallel_4\ \vee \ \parallel \partial_t f \parallel_4\ \leq\ C\ t^{-1/2} \left ( \parallel A_+ \parallel_{4/3}\ + \ \parallel \omega^{-1} \dot{A}_+\parallel_{4/3} \right )
\eeq

\noi while for $\beta$ a multiindex with $\beta = 2$
\beq
\label{2.87e}
\parallel \partial^{\beta} f  \parallel_{4}\ \leq \ C\parallel\Delta f \parallel_4\ = \ 2C \parallel A_0 \parallel_4 \ \leq C \ t^{-1/2}
\eeq

\noi by the Mikhlin theorem, by (\ref{1.10e}) and Lemma 2.4. From (\ref{2.82e}) (\ref{2.83e}) (\ref{2.85e})-(\ref{2.87e}) it follows that  
\beq
\label{2.88e}
\parallel fu_0;L^{8/3}([t, \infty ), W_4^2)\parallel\ \vee \ \parallel \partial_t (fu_0);L^{8/3}([t, \infty ), L^4) \parallel\ \leq\ C\ t^{-13/8} 
\eeq

\noi which together with (\ref{2.84e}) proves that $(fu_0, 0) \in X([1, \infty ))$. \par \nobreak \hfill $\sq$\par

\mysection{The Zakharov system (Z)$_{\bf 2}$ in space dimension n $=$ 2}
\hspace*{\parindent} In this section, we treat the (Z)$_2$ system and
eventually prove Proposition 1.3. As mentioned in the introduction, the
situation is much less satisfactory than in space dimension $n=3$. The
free part $A_0$ of the asymptotic field is estimated at best as
\beq
\label{3.1e}
\parallel A_0(t) \parallel_{\infty} \ \leq C\ t^{-1/2} 
\eeq

\noi and we are unable to handle such a slow decay in Step 1, so that
the final result will eventually be restricted to the special case of
zero asymptotic state $(A_+,\dot{A}_+)$ for $A$. On the other hand, in
a suitable limit, the Zakharov system formally yields the cubic NLS
equation, which is short range for $n=2$, and one might naively expect
a similar situation for the (Z)$_2$ system, allowing for a treatment of
that system without a smallness condition on $u$. This turns out not to
be the case, and the (Z)$_2$ system does actually require such a
smallness condition at the level of Step 1.  The treatment of that step
is very similar to the case of (Z)$_3$. The relevant space $X(\cdot )$
is again given by (\ref{1.6e}), now with $n = 2$, and the main result
can be stated as follows.\\

\noi {\bf Proposition 3.1} {\it Let $h$ be defined as in Section 2 with
$\lambda = 1/2$ and let $X(\cdot )$ be defined by (\ref{1.6e}). Let
$u_a$, $A_a$, $R_1$ and $R_2$ be sufficiently regular and satisfy the estimates
\beq
\label{3.2e}
\parallel u_a(t) \parallel_{\infty} \ \vee \ \parallel \nabla u_a(t) \parallel_{\infty} \ \vee \ \parallel \Delta u_a(t)\parallel_{\infty} \ \vee \ \parallel \partial_t u_a(t) \parallel_{\infty} \ \leq c\ t^{-1} \ ,
\eeq
\beq
\label{3.3e}
\parallel \partial_t^j A_a(t) \parallel_{\infty} \ \leq a \ t^{-1-j\theta} \qquad \hbox{\it for some $\theta >0$ and for $j = 0, 1$}\ ,
\eeq  
\beq
\label{3.4e}
\parallel \partial_t^j R_1; L^1([t , \infty ), L^2)\parallel \ \leq r_1 \ h(t) \qquad \hbox{\it for $j = 0, 1$}\ ,
\eeq 
\beq
\label{3.5e}
\parallel R_1; L^{4}([t , \infty ), L^4)\parallel \ \leq r_1 \ t^{-\eta}\ h(t) \qquad \hbox{\it for some $\eta \geq 0$}\ ,
\eeq 
\beq
\label{3.6e}
\parallel \omega^{-1}R_2; L^1([t , \infty ), H^1)\parallel \ \leq r_2 \ h(t) 
\eeq

\noi for some constants $c$, $a$, $r_1$ and $r_2$ with $c$ sufficiently small and for all $t \geq
1$. Then there exists $T$, $1 \leq T < \infty$, and there exists a
unique solution $(v, B)$ of the system (\ref{1.3e}) in $X([T, \infty
))$.}\\

\noi {\bf Sketch of proof.} The proof is essentially the same as that
of Proposition 2.1 with minor differences in the estimates, and we
concentrate on the latter. We take again $(v, B) \in X([T, \infty ))$
for some $T$, $1 \leq T < \infty$, so that $(v,B)$ satisfies
\bea
\label{3.7e}
&&\parallel v(t) \parallel_2 \ \leq N_0\ h(t)\\
&&\parallel v ; L^{4}(J, L^4)\parallel \ \leq \  N_1\ h(t)\\ 
\label{3.8e}
&&\parallel B(t);H^1 \parallel \ \vee \ \parallel \partial_t B(t)\parallel_2\ \leq \  N_2\ h(t) \\
\label{3.9e}
&&\parallel \partial_t v(t) \parallel_2 \ \leq \  N_3\ h(t)\\
\label{3.10e}
&&\parallel \partial_t v ; L^{4}(J, L^4)\parallel \ \leq \  N_4\ h(t)\\ 
\label{3.11e}
&&\parallel \Delta v(t) \parallel_2\ \leq N_5\ h(t) \\
\label{3.12e}
&&\parallel \Delta v; L^{4}(J, L^4) \parallel\ \leq N_6\ h(t) 
\label{3.13e}
\eea

\noi for some constants $N_i$, $0 \leq i \leq 6$ and for all $t \geq
T$, with $J = [t, \infty )$. We construct a solution $(v',B')$ of the
system (\ref{1.5e}) in $X(I)$ first for $I = [T, t_0]$ and then for $I
= [T, \infty )$. For that purpose we define again
\beq
\label{3.14e}
N'_0 =  \ \mathrel{\mathop {\rm Sup}_{t\in I}}\ h(t)^{-1} \parallel v'(t) \parallel_2
\eeq 
\beq
\label{3.15e}
N'_1 =  \ \mathrel{\mathop {\rm Sup}_{t\in I}}\ h(t)^{-1} \parallel v';L^{4}(J,L^4) \parallel
\eeq
\beq
\label{3.16e}
N'_2 =  \ \mathrel{\mathop {\rm Sup}_{t\in I}}\ h(t)^{-1} \left ( \parallel B'(t);H^1 \parallel \ \vee \ \parallel \partial_t B'(t) \parallel_2 \right )
\eeq
\beq
\label{3.17e}
N'_3 =  \ \mathrel{\mathop {\rm Sup}_{t\in I}}\ h(t)^{-1} \parallel \partial_t v'(t)\parallel_2
\eeq
\beq
\label{3.18e}
N'_4 =  \ \mathrel{\mathop {\rm Sup}_{t\in I}}\ h(t)^{-1} \parallel \partial_t v';L^{4}(J,L^4) \parallel
\eeq
\beq
\label{3.19e}
N'_5 =  \ \mathrel{\mathop {\rm Sup}_{t\in I}}\ h(t)^{-1} \parallel \Delta v'(t)\parallel_2
\eeq
\beq
\label{3.20e}
N'_6 =  \ \mathrel{\mathop {\rm Sup}_{t\in I}}\ h(t)^{-1} \parallel \Delta v';L^{4}(J,L^4) \parallel
\eeq

\noi where $J = [t, \infty ) \cap I$. The crux of the proof is to
estimate the $N'_i$ in terms of the $N_i$. By exactly the same method
as in the proof of Proposition 2.1, we obtain
\beq
\label{3.21e}
N'_0 \leq 2c\ N_2 + r_1
\eeq
\beq
\label{3.22e}
N'_1 \leq C_1 \left ( a \ N'_0 + c \ N_2 + r_1\right )
\eeq

\noi for $N_2 \overline{h}(T) \leq C$,
\beq
\label{3.23e}
N'_3 \leq C_3 \left ( a \ N'_0 \ T^{-\theta} + c\ N_2 + r_1 + \left ( N_2\ N'_4\ N'_1\ \overline{h}(T)\right )^{1/2}\right )
\eeq
\beq
\label{3.24e}
N'_4 \leq C_4 \left ( a \left ( N'_3 + N'_0  \ T^{-\theta}\right )  + c\ N_2 + r_1 +  N_2\left (N'_1 +  N'_4\right )  \overline{h}(T)\right )
\eeq  
\beq
\label{3.25e}
N'_5 \leq 4 \left (  N'_3 + a\ N'_0  \ T^{-1} + c\ N_2 \ T^{-1} + r_1 +  C N'_0 \left (  N_2 h(T)\right )^2\right ) 
\eeq
\beq
\label{3.26e}
N'_6 \leq 4 \left (  N'_4 + a\ N'_1  \ T^{-1} + C c\ N_2 \ T^{-3/4} + r_1 \ T^{-\eta } +  C N'_1 \left (  N_2 h(T)\right )^{4/3}\right ) 
\eeq
\beq
\label{3.27e}
N'_2 \leq C_2 \left ( c \left ( N_0 + N_5  \right )  + N_1 \left (N_1 +  N_6\right )  \overline{h}(T) + r_2\right ) \ .
\eeq 

The estimates (\ref{3.21e})-(\ref{3.27e}) are very similar to the
corresponding estimates (\ref{2.28e}) (\ref{2.30e}) (\ref{2.32e})
(\ref{2.34e}) (\ref{2.36e}) (\ref{2.37e}) (\ref{2.41e}) of the case $n
= 3$. Aside from unimportant changes in the exponents in (\ref{3.25e}) (\ref{3.26e}),
the main differences are (i) the occurrence of the factor $T^{-\theta
}$ in (\ref{3.23e}) and (\ref{3.24e}), coming from the assumption
(\ref{3.3e}) with $j = 1$ on $\parallel \partial_t
A_a\parallel_{\infty}$, and (ii) the replacement of $c T^{-1/2}$ by $c$
everywhere, coming from the assumption (\ref{3.2e}) as compared with
(\ref{2.6e}). The latter difference is responsible for the need of the
smallness condition on $c$. In fact, with the estimates
(\ref{3.21e})-(\ref{3.27e}) available, the proof proceeds as that of
Proposition 2.1. The main step is to prove that the set ${\cal R}$
defined by (\ref{3.7e})-(\ref{3.13e}) is stable under the map $\phi :
(v, B) \to (v', B')$. This is ensured by taking
\beq
\label{3.28e}
\left \{ \begin{array}{l} 
N_0 = 2c\ N_2 + r_1 \\ \\
N_1 = C_1 \left ( a\ N_0 + c\ N_2 + r_1 \right ) \\ \\
N_3 = C_3 \left ( c\ N_2 + r_1 + 1 \right ) \\ \\
N_4 = C_4 \left ( a\ N_3 + c\ N_2 + r_1 + 1 \right )\\ \\
N_5 = 4 \left ( N_3 + r_1 + 1 \right ) \\ \\
N_6 = 4 \left ( N_4 + r_1 + 1 \right )\\ \\
N_2 = C_2 \left ( c \left ( N_0 + N_5 \right ) + r_2 + 1 \right )
\end{array} \right .
\eeq
     
\noi and by taking $T$ sufficiently large so that the remaining $o(1)$
terms in the RHS of (\ref{3.21e})-(\ref{3.27e}) do not exceed 1. In
order to solve the system (\ref{3.28e}), we remark that the constants
$N_1$ and $N_6$ associated with the Strichartz norms do not
occur in the RHS and can therefore be determined at the very end.
Eliminating $N_3$ and $N_4$ (to be determined at the end as functions of $N_2$)
one is left with the reduced system
\beq
\label{3.29e}
\left \{ \begin{array}{l} 
N_0 = 2c\ N_2 + r_1 \\ \\
N_5 = 4 C_3\ c\ N_2 + 4 \left ( C_3 + 1\right ) \left ( r_1 + 1 \right ) \\ \\
N_2 = C_2 \left ( c \left ( N_0 + N_5 \right ) + r_2 + 1 \right )
\end{array} \right .
\eeq

\noi which can obviously be solved for $N_0$, $N_5$ and $N_2$ positive for $c$ sufficiently small.\par

The remaining part of the proof proceeds as that of Proposition 2.1
with appropriate changes in the contraction argument and will be
omitted. \par \nobreak \hfill $\sq$ \par

\noi {\bf Remark 3.1.} The assumption (\ref{3.3e}) on $A_a$ is rather
arbitrary. It is too strong to accomodate a non zero $A_0$ satisfying
only (\ref{3.1e}). On the other hand it is weaker by one power of $t$
than the condition that would be satisfied by an $A_1$ devised to ensure
that $R_2 = 0$. It has been chosen so as to ensure that the proof of
the proposition proceeds smoothly. \\

We are now in a position to derive the final result, namely Proposition
1.3. As already mentioned, the assumption (\ref{3.3e}) forces us to
take $A_0 = 0$. \\

\noi {\bf Proof of Proposition 1.3.}\par

The result will follow from Proposition 3.1 once we have proved that $(u_0, 0)$ satisfies the assumptions of that proposition for $(u_a, A_a)$. From the standard $L^1- L^{\infty}$ estimates of $U(t)$, we obtain
$$\parallel u_0(t) \parallel_{\infty} \ \leq \ (2 \pi t)^{-1} \parallel u_+ \parallel_1$$
$$2 \parallel \partial_t u_0 \parallel_{\infty} \ = \ \parallel \Delta u_0 \parallel_{\infty} \ \leq \ (2 \pi t)^{-1} \parallel \Delta u_+ \parallel_1$$

\noi which proves (\ref{3.2e}). Since $A_a = 0$ and $R_1 = 0$, (\ref{3.3e})-(\ref{3.5e}) are obvious. Finally $R_2 = - \Delta |u_0|^2$, so that (\ref{3.6e}) with $h(t) = t^{-1}$ follows from Lemma 2.3 with $n = 2$.\par \nobreak \hfill $\sq$ \par

\end{document}